\documentclass[journal]{IEEEtran}

\IEEEoverridecommandlockouts                              

\usepackage{multirow}
\usepackage{graphics} 
\usepackage{epsfig} 
\usepackage{subfig}
\usepackage{times} 
\usepackage{amsmath} 
\usepackage{amssymb}  
\usepackage{cite}
\usepackage{color}
\usepackage{algorithm}
\usepackage[noend]{algpseudocode}
\usepackage[switch, pagewise]{lineno}
\renewcommand\linenumberfont{\normalfont\small}

\newtheorem{remark}{Remark}[section]
\newtheorem{lemma}{Lemma}[section]
\newtheorem{theorem}{Theorem}[section]
\newtheorem{corollary}{Corollary}[section]

\graphicspath{{figure/}}

\def\setneig{ \sum_{j \in \mathbf{N}_i } }

\makeatletter
\def\BState{\State\hskip-\ALG@thistlm}
\makeatother

\def\beqn{\begin{eqnarray}}
\def\eeqn{\end{eqnarray}}
\def\pgpd{ p_i^g - p_i^d  }

\title{Distributed coordination for optimal energy generation and distribution in smart grid networks}

\author{\small Hyo-Sung Ahn$^{\dag}$, Byeong-Yeon Kim$^{\ddag}$, Young-Hun Lim$^{\dag}$, and Kwang-Kyo Oh$^{\sharp}$
\thanks{\small $^{\dag}$School of Mechatronics, Gwangju Institute of Science and Technology, Gwangju, Korea.
{E-mails: hyosung@gist.ac.kr, hoonnim@gist.ac.kr}}
\thanks{\small $^{\ddag}$SFR NSSS System Design Division, Korea Atomic Energy Research Institute (KAERI), Daejeon, Korea. {E-mail: byeongyeon@kaeri.re.kr}}
\thanks{\small $^{\sharp}$Automotive Components and Materials R$\&$D Group, Korea Institute of Industrial Technology, Gwangju, Korea. {E-mail: kwangkyo.oh@gmail.com}}
}

\begin{document}
\maketitle


\thispagestyle{empty}
\pagestyle{empty}
\begin{abstract}
This paper proposes coordination laws for optimal energy generation and distribution in energy network, which is composed of physical flow layer and cyber communication layer. The physical energy flows through the physical layer; but all the energies are coordinated to generate and flow by distributed coordination algorithms on the basis of communication information. 
First, distributed energy generation and energy distribution laws are proposed in a decoupled way without considering the interactive characteristics between the energy generation and energy distribution. Then, a joint coordination law is designed in a coupled way taking account of the interactive characteristics. The coordination laws proposed in this paper are fully distributed in the sense that they are decided optimally only using relative information among neighboring nodes. Through numerical simulations, the validity of the proposed distributed coordination laws is illustrated.
\end{abstract}

\section{Introduction}  \label{introduction}
The principal goal of smart grid networks is to supply energy to customers by using distributed energy resources (DERs), in order to reduce greenhouse gas emissions and to save cost arising when generating energy in a centralized way \cite{Tuballa_2016_RSER}. 
The smart grid networks may be called cyber physical energy systems \cite{Andres_2009_ISGT}. The smart grid networks are composed of physical energy resources and communication networks. The energy flows through physical network, while information is exchanged through the cyber communication network. The energy flow through physical network would take a much longer time, while the information through the cyber layer is exchanged extremely fast by the benefit of the advance of computational power and ultra high channel capacity in communication network. In renewable energy networks, it is a fundamental but tough issue to make a balance between power supply and demand \cite{Spataru_2012_ISGT}, particularly due to the intermittent nature of distributed energy resources (DERs) \cite{Kumar_2015_TSG}. In smart grid, each node acts as a customer and generator; so each node has a desired energy level while nodes have different generation cost. Thus, the overall network attempts to achieve the supply-demand balance by generating the total energy required in the network, while taking account of the different generation and distribution costs. 


As energy exchange becomes active making use of distributed energy resources, optimization in the direction of saving costs in generation and distribution becomes more and more important. 
However, due to the lack of information in cost estimation in a global sense and economic benefits for global optimal generation, it is difficult to achieve optimization in a centralized way \cite{Hidalgo_2011_ISGT}. That is, although renewable energy has become popular for replacing
traditional energy resources, the
uncertainty associated with the operation has
become a challenge \cite{Harirchi_2014_ICSGC}, which implies that the autonomous generation with only local and distributed information is highly required. Note that in most of existing works, for coordinated
control of the distributed generators integrated to a DC microgrid, feedback and/or feed-forward control algorithms have been developed in dynamic levels \cite{Kumar_2015_TSG,Xin_2013_CTA}. 
In microgrid or in smart grid, decentralized control for energy distribution has become feasible due to multiterminal DC systems \cite{Gavriluta_2014_tia}. The energy distribution concept may be connected to the energy trading, which may be conducted by coordinating energy flow such as using vehicle-to-grid (V2G) services \cite{Al-Awami_2012_TSG}. 
So, traditional distributed optimization algorithms \cite{Johansson_2008_KTH,Nedic_etal_2010_TAC} may be used to find optimal path in a distributed way; but, the existing distributed optimization algorithms are not specified to the energy generation and distribution. Thus, it may take a huge amount of computation, when using nonlinear distributed optimization schemes. As alternative solutions, distributed control and/or consensus approach may be utilized for energy distribution 
\cite{Lee_Mesbahi_2011_CDC,Kar_Hug_2012_PESGM,Tychogiorgos_2013_TWC}; but existing distributed approaches are mostly specified to state re-assignment problems without considering supply-demand balance from a perspective of cyber physical systems. 

The main purpose of this paper is to propose energy generation and distribution algorithms in cyber physical energy networks under the constraint of meeting the supply-demand balance in a consensus approach. Note that this problem has been already studied in \cite{Kim_2015_CEP}; but in \cite{Kim_2015_CEP}, they have not considered optimization issue. 
As related works, in
\cite{Kim_2015_TCST}, consensus-based
coordination and control for temperature control was proposed; but in \cite{Kim_2015_TCST}, the interactive characteristics between energy generation and energy flow was not investigated. In \cite{WuYuan_2014_CCC}, they designed a schedule algorithm for energy generation in distributed energy networks to minimize the energy provisioning
cost. 
However, in \cite{WuYuan_2014_CCC}, there is no consideration of supply-demand balance, and the optimization is solved in a centralized way. In \cite{Nottrott_2012_PESGM}, an optimal energy management system was proposed for dispatch scheduling to avoid peak loads while considering storage capacity. 

Contrary to the existing works, 
as the contribution of this paper, we first design an optimal generation law in a simple manner using algebraic characteristics of Laplacian matrix, while meeting the supply-demand balance under a distributed coordination setup. Second, optimal distribution algorithm is designed taking account of the cost in energy flows. The costs are assigned to edges in network graph. The optimal generation and optimal distribution laws are realized in a fully distributed way. 
Third, we propose a joint optimal coordination law taking account of both generation cost and distribution cost in a coupled manner. Thus, by realizing joint coordination law, we can minimize the overall cost while considering both the generation and energy flow simultaneously.

\section{Problem Formulation} \label{formulation}
\begin{figure}
\includegraphics[width=0.5\textwidth]{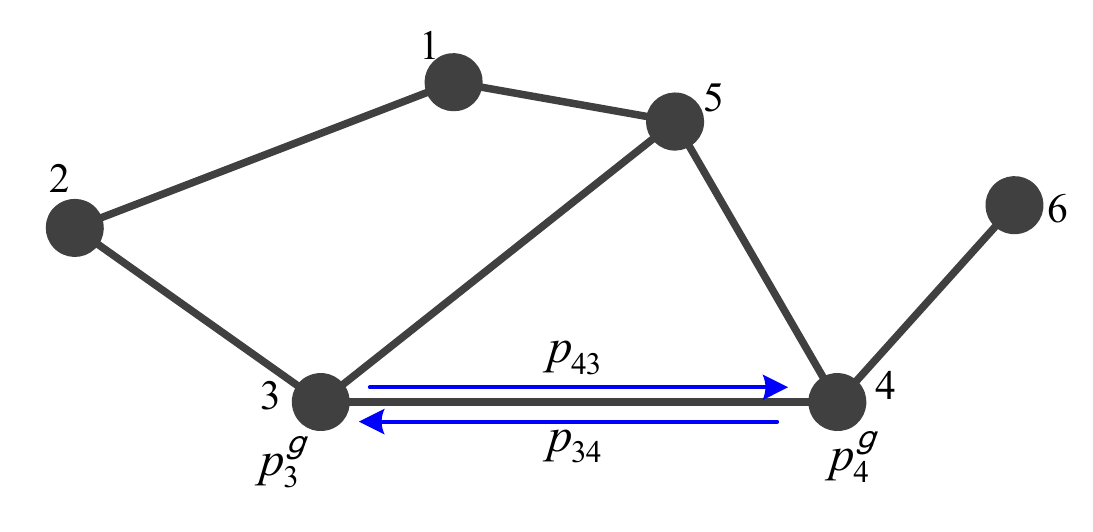}
\caption{Energy network: Nodes are represented by dark circles and interconnections are represented by solid lines.}
\label{sgn}
\end{figure}

The smart grid network is composed of a number of nodes and the nodes are interconnected by physical and cyber layer networks. 
The nodes can be considered as vertices while the interconnections could be considered as edges from the concepts of graph theory (see Fig.~\ref{sgn}). So, we would represent the smart grid network by a graph, which is a triple as $\mathcal{G} = (\mathcal{V},\mathcal{E},\mathcal{W})$, with cardinalities of $\vert \mathcal{V} \vert =n$ and $\vert \mathcal{E} \vert =m$, where $\mathcal{W}$ is the set of weights assigned to each edges.  The cost arising when passing through the physical lines in network would be defined by weights. 
Energy generation in each node is represented by $p_i^g$, and energy flow from the $j$-th node to the $i$-th node is represented by $p_{ij}$; so it is natural to have the constraints of $\vert p_{ij} \vert = \vert p_{ji} \vert$ and $\text{sign}(p_{ij}) =-\text{sign}(p_{ji})$. The energy level of each node is defined as 
\begin{align}
p_i = p_i^g + \sum_{j \in {\bf N}_i} p_{ij} \label{basic_eq}
\end{align} 
where $p_i^g$ is the generated energy at $i$-th node and $\sum_{j \in {\bf N}_i} p_{ij}$, where ${\bf N}_i$ is the set of neighboring nodes of the $i$-th node, is net energy sum exchanged with neighboring nodes. Given the desired energy of each node as $p_i^d$, it should be true that $\sum_{i=1}^n p_i^g = \sum_{i=1}^n p_i^d$, which is called supply-demand balance in smart grid networks as mentioned in Section~\ref{introduction}.  

If each node has enough ability for energy generation or if the cost of energy generation of each node is cheap enough, then given desired energy level $p_i^d$ in each node, each node may simply generate $p_i^g$ such that it would be equal to $p_i^d$. However, if the cost of energy generation of the nodes is different, depending on usability of resource or environment, each node may not be able to generate the desired energy. At some nodes, the generation cost may be expensive while in some other nodes, the cost may be cheap. So, the overall task of distributed generation is to generate $p_i^g$ such that $\sum_{i=1}^n p_i^g=\sum_{i=1}^n p_i^d$ with a minimum cost in the network. Thus, the generated energy in each node may be not equal to the desired energy level, i.e., $p_i^g \neq p_i^d$. The cost of energy generation of each node is given in a convex function form as $\xi_i (p_i^g)^2 + \zeta_i p_i^g + \eta_i$, where the parameters $\xi_i, \zeta_i$, and $\eta_i$ are only available to the $i$-th node. So, we now have the following optimization problem:
\begin{eqnarray}
&\text{min} \left\{ \sum_{i \in \mathcal{V}} (\xi_i (p_i^g)^2 + \zeta_i p_i^g + \eta_i) \right\},~\forall i \in \mathcal{V}& \label{main1_costfunc}\\
&\text{s.t.}& \nonumber\\
& \sum_{i=1}^n p_i^g=\sum_{i=1}^n p_i^d & \label{generation_constraint1}
\end{eqnarray}
Using the following Lagrangian 
\begin{align}
\mathfrak{L}^g &=   \sum_{i \in \mathcal{V}} (\xi_i (p_i^g)^2 + \zeta_i p_i^g + \eta_i) + \lambda (\sum_{i=1}^n p_i^g - \sum_{i=1}^n p_i^d) 
\end{align}
we can have two constraints 
\begin{eqnarray}
&2\xi_i p_i^g + \zeta_i + \lambda =0& \label{g_const1}\\
&\sum_{i=1}^n p_i^g =\sum_{i=1}^n p_i^d& \label{g_const2}
\end{eqnarray}
Thus, for a distributed coordination for optimal energy generation, we need to solve the above equations in a distributed way. 

Meanwhile, the goal of flow coordination for energy distribution is to determine energy flow $p_{ij}$ between neighboring nodes such as $p_i \rightarrow p_i^d$ since $p_i^g \neq p_i^d$. When energy is flowing between neighbor nodes, it is supposed that there are some cost, which can be a function of distance or network capacity; so, it would be desirable to make energy flow more along the edges that take less cost. Eventually, the cost of a  path $(i,j)$ could be defined as a function of $p_{ij}$ such as $g(p_{ij}) \triangleq \alpha_{ij} p_{ij}^2 + \beta_{ij} p_{ij} + \gamma_{ij}$ for the edge $(i,j)$. Here, it is supposed that $\alpha_{ij} = \alpha_{ji}$, $\beta_{ij} = -\beta_{ji}$, and $\gamma_{ij} = \gamma_{ji}$, where the parameters $\alpha_{ij}, \beta_{ij}$, and $\gamma_{ij}$ are available to the pair of the neighboring nodes $i$ and $j$. Then, the overall energy flow under the direction of minimizing the cost can be formulated as
\begin{eqnarray}
&\text{min} \left\{ \sum_{(i,j) \in \mathcal{E}} (\alpha_{ij} p_{ij}^2 + \beta_{ij} p_{ij} + \gamma_{ij}) \right\},~\forall (i,j) \in \mathcal{E}& \label{main_costfunc}\\
&\text{s.t.}& \nonumber\\
& p_i \Rightarrow p_i^d = p_i^g + \sum_{j\in \mathbf{N}_i} p_{ij}, ~\forall i \in \mathcal{V} & \label{flow_constraint1}\\
& p_{ij} + p_{ji} =0,~\forall (i,j) \in \overline{ \mathcal{E} } & \label{flow_constraint2}
\end{eqnarray}
where $\mathcal{E}= \mathcal{E}^+ \cup \mathcal{E}^-$, and $\overline{\mathcal{E}}$ could be $\mathcal{E}^+$ or $\mathcal{E}^-$; here $\mathcal{E}$ is the set of edges of undirected graph, while $\mathcal{E}^+$ is the set of edges when the edges of the same graph are assigned directions arbitrarily; so $\mathcal{E}^-$ is the set of edges with opposite directions. To solve the above optimization, let us use the following Lagrangian:
\begin{align}
\mathfrak{L}^e &=  \sum_{i}^n \sum_{j \in \mathbf{N}_i} g(p_{ij}) + \sum_{i=1}^n \left\{   \lambda_i^v  \left[ (\sum_{j \in \mathbf{N}_i} p_{ij}) + p_i^g -p_i^d  \right] \right\} \nonumber\\
&~~~~~+ \sum_{(i,j) \in \overline{ \mathcal{E}} } \left\{ \lambda_{(i,j)}^{e} (p_{ij} + p_{ji}) \right\}
\end{align}
By conducting derivatives of the Lagrangian $\mathfrak{L}^e$ with respect to $p_{ij}$, $\lambda_i^v$, and $\lambda_{(i,j)}^{e}$, we obtain three main constraints 
\begin{eqnarray}
& 2 \alpha_{ij} p_{ij} + \beta_{ij} + \lambda_i^v + \lambda_{(i,j)}^{e} =0, ~\forall (i,j) \in \mathcal{E}& \label{lagrangian_pij}\\
& \setneig p_{ij} + p_i^g - p_i^d = 0,~\forall i \in \mathcal{V}& \label{lagrangian_lamv}\\
& p_{ij} + p_{ji} = 0,~\forall (i,j) \in \overline{\mathcal{E}} &  \label{lagrangian_lame}
\end{eqnarray}
From (\ref{lagrangian_pij}), obtaining $p_{ij}= -\frac{1}{2 \alpha_{ij} } [\beta_{ij} + \lambda_i^v + \lambda_{(i,j)}^{e}]$, we insert it into (\ref{lagrangian_lamv}) and (\ref{lagrangian_lame}) respectively to have 
\begin{eqnarray}
&-\setneig \left\{   \frac{1}{2 \alpha_{ij} } [\beta_{ij} + \lambda_i^v + \lambda_{(i,j)}^{e}]  \right\} + \pgpd = 0 & \label{problem_1}\\
& \frac{1}{2 \alpha_{ij} } [\beta_{ij} + \lambda_i^v + \lambda_{(i,j)}^{e}] + \frac{1}{2 \alpha_{ji} } [\beta_{ji} + \lambda_j^v + \lambda_{(j,i)}^{e}] =0 \label{problem_2_pre}& 
\end{eqnarray}
So, for a distributed coordination for optimal energy flow, we need to solve the above equations relying upon only neighboring interactions. It is observed that the optimization problem defined in (\ref{main1_costfunc})-(\ref{generation_constraint1}), and the optimization problem defined in (\ref{main_costfunc})-(\ref{flow_constraint2}) can be solved independently. However, since in (\ref{flow_constraint1}), it is also a function of $p_i^g$, depending upon the amount of $p_i^g$ at each node,  the optimization cost of (\ref{main_costfunc}) might be changing. In the next section, we solve these optimization problems both without considering the interactive characteristics (i.e., we first solve the optimization (\ref{main1_costfunc})-(\ref{generation_constraint1}); then, with the generated $p_i^g$, we solve (\ref{main_costfunc})-(\ref{flow_constraint2})) and with a consideration of the interactive characteristics (i.e., we solve the two optimization problems simultaneously).

\section{Distributed Solutions}\label{main_laplacian}
In this section, we solve the optimization problems for energy generation and energy distribution defined in Section~\ref{formulation} in a distributed way.  
\subsection{Coordination of energy generation} \label{subsec_gen}
Let us insert $p_i^g$ in (\ref{g_const1}) into (\ref{g_const2}) to obtain:
\begin{align}
\lambda= -\frac{  \sum_{i=1}^n \left[ p_i^d + \frac{\zeta_i}{2 \xi_i}       \right]   }{\sum_{i=1}^n \frac{1}{2\xi_i}} \label{g_const3}
\end{align}
The Lagrangian parameter $\lambda$ in the above equation can be obtained if there is a central computational unit since all the information could be gathered. To obtain $\lambda$ in a distributed way, let us introduce new subsidiary variables $v_i(t)$ and $u_i(t)$ and assign their initial values as 
\begin{align}
v_i(0) &= \frac{1}{2 \xi_i} \label{vi_initial_assign} \\
u_i(0) &= p_i^d + \frac{\zeta_i}{2 \xi_i} \label{ui_initial_assign}
\end{align}
Let us update $v_i$ and $u_i$ by the following distributed consensus algorithms:
\begin{align}
\dot{v}_i &= \sum_{j \in \mathbf{N}_i} (v_j - v_i) \label{consensus_vi} \\
\dot{u}_i &= \sum_{j \in \mathbf{N}_i} (u_j - u_i) \label{consensus_ui} 
\end{align}
For a distributed coordination, we use the following average consensus law \cite{olfati2007consensuspieee}:
\begin{lemma}  \label{average_consensus}
If the network defined by the graph $\mathcal{G}$ is connected, the above consensus laws ensure the convergence of states to the following average value:
\begin{align}
v_1 = v_2 = \cdots = v_n \rightarrow  \frac{1}{n} \sum_{i=1}^n (v_i(0))  \nonumber \\
u_1 = u_2 = \cdots = u_n \rightarrow  \frac{1}{n} \sum_{i=1}^n (u_i(0))   \nonumber
\end{align}
\end{lemma}

Based on the above lemma, through the paper, it is supposed that the network is connected. 
Now, we can make the following result.
\begin{theorem} \label{theorem_energy_generation_main}
Let the converged values of (\ref{consensus_vi}) and (\ref{consensus_ui}) be denoted as $v_i^\ast$ and $u_i^\ast$, respectively. Then, $\lambda^\ast$ is given as
\begin{align}
\lambda^\ast = -\frac{u_i^\ast}{v_i^\ast} \label{main1_lambda}
\end{align}
and it is equal to (\ref{g_const3}). 
\end{theorem}
\begin{IEEEproof}
Since (\ref{consensus_vi}) and (\ref{consensus_ui})  are average consensus, it is clear by the \textit{Lemma~\ref{average_consensus}} that $v_i^\ast = \frac{1}{n}\left[  \sum_{i=1}^n \frac{1}{2 \xi_i} \right]$ and $u_i^\ast = \frac{1}{n}\left[  \sum_{i=1}^n (p_i^d + \frac{\zeta_i}{2 \xi_i}) \right]$. Thus, inserting $v_i^\ast = \frac{1}{n}\left[  \sum_{i=1}^n \frac{1}{2 \xi_i} \right]$ and $u_i^\ast = \frac{1}{n}\left[  \sum_{i=1}^n (p_i^d + \frac{\zeta_i}{2 \xi_i}) \right]$ into (\ref{main1_lambda}), we see that $\lambda^\ast$ of  (\ref{main1_lambda}) becomes (\ref{g_const3}). 
\end{IEEEproof}

The Lagrangian parameter $\lambda^\ast$ in (\ref{main1_lambda}) could be obtained in a distributed way; so after obtaining it, when each node generates energy such as
\begin{align}
p_i^g = -\frac{1}{2\xi_i}(\lambda^\ast + \zeta_i) \label{p_i_g_final}
\end{align}
the cost function (\ref{main1_costfunc}) will be minimum while the equality constraint (\ref{generation_constraint1}) would be satisfied. The overall process for energy generation is summarized in the \textit{Algorithm~\ref{coordination1}}.
\begin{algorithm}
\caption{Coordination for generation}\label{coordination1}
\begin{algorithmic}[1]
\State $\textit{assign initial values of}$ $v_i(0)$ \textit{and} $u_i(0)$ using (\ref{vi_initial_assign}) and  (\ref{ui_initial_assign})
\State $\textit{conduct distributed consensus}$ for $v_i$ and $u_i$ using (\ref{consensus_vi}) and (\ref{consensus_ui})
\State with converged $\lambda^\ast$, \textit{compute} $p_i^g$ 
using (\ref{p_i_g_final})
\end{algorithmic}
\end{algorithm}

In some cases, it may be strictly required to set upper and lower boundaries in the amount of energy generation of each node. That is, we may additionally have constraints of $\underline{p_i^g} \leq p_i^g \leq \overline{p_i^g}$. In what follows, we provide an analysis associated with this issue. By inserting $\lambda^\ast$ in (\ref{main1_lambda}) into (\ref{p_i_g_final}), we have
\begin{align}
p_i^g = \frac{1}{2 \xi_i} \left\{  \frac{ \sum_{i=1}^n p_i^d     }{ \sum_{i=1}^n \frac{1}{2 \xi_i}  }  + \frac{\sum_{i=1}^n \frac{\zeta_i}{2 \xi_i}  }{   \sum_{i=1}^n \frac{1}{2 \xi_i}     }   - \zeta_i    \right\} \label{pigg}
\end{align}
From the above equation, as some special instances, if $\xi_1 = \xi_2 = \cdots = \xi_n$ and $\zeta_1 = \zeta_2 = \cdots = \zeta_n$, then $p_i^g$ is the average of the desired energy, i.e., $p_i^g = \frac{1}{n} \sum_{i=1}^n p_i^d$. Further, denoting $p^D \triangleq \sum_{i=1}^n p_i^d$, if only $\zeta_1 = \zeta_2 = \cdots = \zeta_n$ holds, then $p_i^g$ is obtained as $p_i^g  =  \frac{p^D}{\sum_{j=1}^n \frac{\xi_i}{\xi_j} }$. 
Thus, if $\xi_i, i=1, \cdots, n$ are selected as
\begin{align}
\frac{p^D}{\overline{p_i^g}}  \leq  \sum_{j=1}^n \frac{\xi_i}{\xi_j} \leq \frac{p^D}{\underline{p_i^g}}
\end{align}
then the lower and upper boundary conditions $\underline{p_i^g} \leq p_i^g \leq \overline{p_i^g}$ will be satisfied. In general cases, i.e., $\xi_i \neq \xi_j$ and  $\zeta_i \neq \zeta_j$, then, with time-variant variables $x_i(t), y_i(t), z_i(t)$, we set their initial values as
$x_i(0) = p_i^d$, $y_i(0) = \frac{1}{2\xi_i}$, and $z_i(0) = \frac{\zeta_i}{2 \xi_i}$.
Then, we update $x_i, y_i, z_i$ by $\dot{x}_i = \sum_{j \in \mathbf{N}_i} (x_i - x_j)$, 
$\dot{y}_i = \sum_{j \in \mathbf{N}_i} (y_i - y_j)$, and $\dot{z}_i = \sum_{j \in \mathbf{N}_i} (z_i - z_j)$.
Denoting the converged solutions as $x_i^\ast, y_i^\ast, z_i^\ast$, we can compute the right-hand side of (\ref{pigg}) as ${p_i^g}^\ast = \frac{1}{2\xi_i} (x_i^\ast/y_i^\ast + z_i^\ast/y_i^\ast - \zeta_i)$. Thus, if 
$\underline{p_i^g} \leq {p_i^g}^\ast \leq \overline{p_i^g}$, it can be supposed that the energy will have been generated within the lower and upper boundaries.

\subsection{Coordination of energy flow} \label{subsec_flow}
Let us denote $w_i \triangleq -\setneig \left\{   \frac{1}{2 \alpha_{ij} } \beta_{ij}   \right\} +(\pgpd)$. Since $\alpha_{ij} = \alpha_{ji}$ and $\beta_{ij}=-\beta_{ji}$, we can rewrite the equations (\ref{problem_1}) and (\ref{problem_2_pre}) as 
\begin{eqnarray}
& \setneig    \frac{1}{2 \alpha_{ij} } (\lambda_i^v + \lambda_{(i,j)}^{e}) = w_i & \label{vertex_lambda}\\
& \lambda_i^v +\lambda_j^v+ 2 \lambda_{(i,j)}^{e} = 0& \label{edge_lambda}
\end{eqnarray}
Now, inserting  $\lambda_{(i,j)}^{e}$ from  (\ref{edge_lambda}) into (\ref{vertex_lambda}), we have 
\begin{eqnarray}
 \setneig  \frac{1}{4 \alpha_{ij}}  (\lambda_i^v -\lambda_j^v) = w_i
\end{eqnarray}
which is equivalent to the following vector form
\begin{eqnarray}
\mathcal{L}_{\mathcal{G}} \lambda^v = {\bf w} \label{laplacian_lambda_v}
\end{eqnarray}
where $\mathcal{L}_{\mathcal{G}}$ is the weighted Laplacian of the graph, $\lambda^v=[\lambda_1^v, \cdots, \lambda_n^v]^T$, and ${\bf w}=[-\sum_{j \in \mathbf{N}_1} \left\{   \frac{1}{2 \alpha_{1j} } \beta_{1j}   \right\} +(p_1^g- p_1^d), \cdots, -\sum_{j \in \mathbf{N}_n} \left\{   \frac{1}{2 \alpha_{nj} } \beta_{nj}   \right\} +(p_n^g- p_n^d)]^T$. The weighted Laplacian $\mathcal{L}_{\mathcal{G}}$ is computed as $\mathcal{L}_{\mathcal{G}} = D_{\mathcal{G}} - A_{\mathcal{G}}$ where adjacency matrix is calculated as
\begin{align}
A_{\mathcal{G}}(i,j) = \begin{cases} \frac{1}{4 \alpha_{ij}} ~\text{if}~(i,j) \in \mathcal{E} \\ 0 ~\text{otherwise}    \end{cases}
\end{align}
and the degree matrix is given as $D_{\mathcal{G}}(i,i) = \setneig A_{\mathcal{G}}(i,j)$. 
Also, from (\ref{edge_lambda}), using the incidence matrix $\mathcal{H}$, we can obtain the following relationship:
\begin{eqnarray}
\lambda^e = -\frac{1}{2}\mathcal{H}^T \lambda^v  \label{laplacian_lambda_e}
\end{eqnarray}
where $\lambda^e=[\cdots,  \lambda_{(i,j)}^{e}, \cdots ]^T$.
It is observed that, from (\ref{laplacian_lambda_v}) and (\ref{laplacian_lambda_e}), if we could solve (\ref{laplacian_lambda_v}) in a distributed way, then $\lambda^{e}$ can be obtained in a distributed way directly. Thus, the key issue now becomes to solving (\ref{laplacian_lambda_v}) in a distributed way. It is well-known that the Laplacian matrix has simple zero eigenvalue as its one of the eigenvalues; so it is singular. Thus, (\ref{laplacian_lambda_v}) cannot be solved by taking an inverse for any case. 
To proceed further, let us define an orthogonal matrix $\mathcal{Q}\in \Bbb{R}^{n\times n}$ such that
\begin{eqnarray}
\mathcal{Q} =\left[ \begin{matrix}  q_1 & q_2 &\cdots & q_n  \end{matrix} \right] \label{Trans_Q}
\end{eqnarray}
where $q_i \in \Bbb{R}^{n \times 1}$ and $q_i^T q_j =0$ and $q_1 = \epsilon {\bf 1}_n$, where $\epsilon >0$ is a positive constant and ${\bf 1}_n$ is the length $n$ vector with all elements being one. We now define a new Lagrangian vector $\overline{\lambda}^v$ as
\begin{align}
{\lambda}^v  = \mathcal{Q} \overline{\lambda}^v \label{new_lambda_trans_Q}
\end{align}
For the main result, from (\ref{laplacian_lambda_v}), let us propagate $\lambda^v$ such as
\begin{align}
\dot{\lambda}^v = -\mathcal{L}_{\mathcal{G}} \lambda^v + {\bf w} \label{laplacian_lambda_v_prop}
\end{align}
\begin{lemma} \label{lemma_const}
There exists a constant $\lambda^v$ such that (\ref{laplacian_lambda_v}) holds if and only if $\dot{\lambda}^v \rightarrow 0$ in (\ref{laplacian_lambda_v_prop}).
\end{lemma}
\begin{IEEEproof}
It is clear that if $\dot{\lambda}^v \rightarrow 0$, then $\mathcal{L}_{\mathcal{G}} \lambda^v = {\bf w}$. Let us suppose that $\dot{\lambda}^v \neq 0$, i.e., $\dot{\lambda}^v = \psi$; then clearly, $\lambda^v$ is not constant. 
\end{IEEEproof}

\begin{theorem} The solution of (\ref{laplacian_lambda_v_prop}) converges to a constant vector, and the converged solution is equivalent to the solution of $\mathcal{L}_{\mathcal{G}} \lambda^v - {\bf w} =0$. \label{theorem_flow}
\end{theorem}
\begin{IEEEproof}
Let us transform (\ref{laplacian_lambda_v_prop}) using $\overline{\lambda}^v$ as
\begin{align}
\dot{\overline{\lambda}^v} = - \mathcal{Q}^T \mathcal{L}_{\mathcal{G}} \mathcal{Q} \overline{\lambda}^v + \mathcal{Q}^T {\bf w} \label{laplacian_lambda_v_prop2}
\end{align}
Using the fact that $({\bf 1}_n)^T {\bf w}=0$ and using the property of Laplacian matrix, we can change (\ref{laplacian_lambda_v_prop2}) as
\begin{equation}
\dot{\overline{\lambda}^v} = -\left[
\begin{array}{c|c}
  0 & 0 \\
  \hline
  0 & \mathcal{Q'}^T \mathcal{L}_{\mathcal{G}} \mathcal{Q'}
\end{array}
\right]  \overline{\lambda}^v   +    \left[
\begin{array}{c}
  0  \\
  \hline
  {\bf w}' 
\end{array}
\right]        \label{laplacian_lambda_v_prop3}
\end{equation}
where $\mathcal{Q}' \in \Bbb{R}^{n \times (n-1)}$ is the reduced orthogonal matrix given as
\begin{eqnarray}
\mathcal{Q}' =\left[ \begin{matrix}  q_2 &\cdots & q_n  \end{matrix} \right] \label{q_prime}
\end{eqnarray}
and ${\bf w}'$ is the length $n-1$ vector determined by ${\bf w}'= (\mathcal{Q}')^T {\bf w}$. It is noticeable that all eigenvalues of $\mathcal{Q'}^T \mathcal{L}_{\mathcal{G}} \mathcal{Q'}$ are positive real, i.e., it is positive definite. Thus, it is obvious that from (\ref{laplacian_lambda_v_prop3}), ${\overline{\lambda}^v}$ is convergent. That is, ${\overline{\lambda}^v} \rightarrow {\overline{\lambda}^v}^\ast$ as $t$ increases. Furthermore, from ${\lambda}^v  = \mathcal{Q} \overline{\lambda}^v$, we have ${{\lambda}^v}^\ast  = \mathcal{Q} {\overline{\lambda}^v}^\ast$.
Consequently, by the \textit{Lemma~\ref{lemma_const}}, the solution of (\ref{laplacian_lambda_v_prop}) is equivalent to the solution of $\mathcal{L}_{\mathcal{G}} \lambda^v - {\bf w} =0$. 
\end{IEEEproof}

The propagation of (\ref{laplacian_lambda_v_prop}) can be done in a distributed way such as:
\begin{align}
\dot{\lambda}_i^v = -\setneig \frac{1}{4\alpha_{ij}} (\lambda_i^v - \lambda_j^v) -\setneig \left\{   \frac{1}{2 \alpha_{ij} } \beta_{ij}   \right\} +(p_i^g - p_i^d) \label{final_lambda_v}
\end{align}
Denote the solution as ${{\lambda}_i^v}^\ast$; then, the Lagrangian parameters of edges can be obtained as:
\begin{align}
{\lambda_{(i,j)}^{e ^\ast}} = -\frac{1}{2}( {{\lambda}_i^v}^\ast + {{\lambda}_j^v}^\ast) \label{final_lambda_edge}
\end{align}
Finally, we can make the following main result:
\begin{theorem}
The energy flow 
\begin{align}
 p_{ij} &= -\frac{1}{2 \alpha_{ij}}( \beta_{ij} + {\lambda_i^v}^\ast + {\lambda_{(i,j)}^{e ^\ast}}) \nonumber\\&=-\frac{1}{2 \alpha_{ij}}( \beta_{ij} + \frac{1}{2}{\lambda_i^v}^\ast - \frac{1}{2}{\lambda_j^v}^\ast )  \label{energyflowij}
\end{align}
achieves the minimum cost in the optimization problem (\ref{main_costfunc}). 
\end{theorem}
\begin{IEEEproof}
Inserting $p_{ij}$ from (\ref{energyflowij}) into (\ref{lagrangian_pij}) and (\ref{lagrangian_lame}), we can see that (\ref{lagrangian_pij}) and (\ref{lagrangian_lame}) are satisfied. Also from (\ref{final_lambda_v}), since $\dot{\lambda}_i^v =0$, we have
\begin{align}
 -\setneig \frac{1}{4\alpha_{ij}} (\lambda_i^v - \lambda_j^v) -\setneig \left\{   \frac{1}{2 \alpha_{ij} } \beta_{ij}   \right\} +(p_i^g - p_i^d) =0 \label{final_lambda_edge2}
\end{align}
Then, also by inserting both $p_{ij}$ from (\ref{energyflowij}) and (\ref{final_lambda_edge2}) into (\ref{lagrangian_lamv}), we can see that  (\ref{lagrangian_lamv}) is satisfied. Consequently, since the cost function $g(p_{ij})$ is of convex, we achieve the global minimum of the optimization problem (\ref{main_costfunc}).  
\end{IEEEproof}

The \textit{Algorithm~\ref{coordination2}} summarizes the procedure for energy distribution.
\begin{algorithm}
\caption{Coordination for energy distribution}\label{coordination2}
\begin{algorithmic}[1]
\BState \emph{Given generated energy and desired energy level}, $p_i^g$ and $p_i^d$,
\State $\textit{propagate}$ $\lambda_i^v$ according to (\ref{final_lambda_v}) in a distributed way
\State \textit{exchange} the converged solution ${{\lambda}_i^v}^\ast$ with neighbor nodes
\State \textit{compute} $\lambda_{(i,j)}^{e ^\ast}$ using (\ref{final_lambda_edge})
%
%
\State \textit{compute}  $p_{ij}$ using  ${\lambda_i^v}^\ast$ and  ${\lambda_{(i,j)}^{e ^\ast}}$ by      (\ref{energyflowij}) 
\end{algorithmic}
\end{algorithm}

\subsection{Joint coordination for energy generation and flow} \label{subsec_gen_flow}
Here, we would like to decide $p_i^g$ and $p_{ij}$ in an optimal way under the direction of minimizing (\ref{main1_costfunc}) and (\ref{main_costfunc}) simultaneously. Since these two cost functions can be combined into a convex optimization, with the constraints of (\ref{generation_constraint1}), (\ref{flow_constraint1}), and (\ref{flow_constraint2}), we can define the following Lagrangian:
\begin{align}
\mathfrak{L} &=  \sum_{i \in \mathcal{V}} (\xi_i (p_i^g)^2 + \zeta_i p_i^g + \eta_i) \nonumber\\&~+\sum_{i}^n \left\{ \sum_{j \in \mathbf{N}_i} [ \alpha_{ij} p_{ij}^2 + \beta_{ij} p_{ij} + \gamma_{ij}    ] \right\}  \nonumber\\&~ +  \lambda (\sum_{i=1}^n p_i^g - \sum_{i=1}^n p_i^d)  \nonumber\\&~ + \sum_{i=1}^n \left\{   \lambda_i^v  \left[ (\sum_{j \in \mathbf{N}_i} p_{ij}) + p_i^g -p_i^d  \right] \right\} \nonumber\\&~
+ \sum_{(i,j) \in \overline{ \mathcal{E}} } \left\{ \lambda_{(i,j)}^{e} (p_{ij} + p_{ji}) \right\}
\end{align}
By taking gradients with respect to $p_i^g, p_{ij}, \lambda,   \lambda_i^v$, and $\lambda_i^e$, we have the four same constraints as (\ref{g_const2}), (\ref{lagrangian_pij}), (\ref{lagrangian_lamv}), (\ref{lagrangian_lame}), and the following additional constraint
\begin{align}
2 \xi_i p_i^g + \zeta_i + \lambda + \lambda_i^v =0, ~\forall i \in \mathcal{V}, \label{additional}
\end{align}
Then, from (\ref{additional}), we insert $p_i^g = -\frac{1}{2 \xi_i}[ \zeta_i + \lambda + \lambda_i^v]$ into (\ref{g_const2}) to have
\begin{align}
\lambda = \frac{  \sum_{i=1}^n  \left\{   \frac{\zeta_i}{2 \xi_i} + \frac{ \lambda_i^v }{2\xi_i}   + p_i^d             \right\}   }{  -\sum_{i=1}^n  \frac{1}{2\xi_i}   } \label{lambda_optimal}
\end{align}
In the joint coordination, the main difficulty arises from the fact that $p_i^g$ in (\ref{lagrangian_lamv}) also needs to be decided. So, we cannot solve (\ref{lagrangian_pij}), (\ref{lagrangian_lamv}), (\ref{lagrangian_lame}) to find an optimal $\lambda_i^v$ using the process given in Subsection~\ref{subsec_flow}. 

It is however certain that, given $\lambda_i^v$, we can compute $\lambda$ in (\ref{lambda_optimal}) in a distributed way. The overall idea is same as the approach in Subsection~\ref{subsec_gen}. To find optimal $\lambda$  in a distributed way, we use subsidiary variables $r_i(t)$ and $s_i(t)$ and assign their initial values as $r_i(0) = \frac{1}{2 \xi_i}$ and $s_i(0) =   \frac{\zeta_i}{2 \xi_i} +   \frac{ {{\lambda}_i^v}^\ast }{ 2 \xi_i} +  p_i^d$ respectively. Let us update $r_i$ and $s_i$ by the following distributed consensus algorithms:
\begin{align}
\dot{r}_i &= \sum_{j \in \mathbf{N}_i} (r_i - r_j) \label{consensus_yi} \\
\dot{s}_i &= \sum_{j \in \mathbf{N}_i} (s_i - s_j) \label{consensus_zi} 
\end{align}
Now, we can make the following result.
\begin{corollary}
Let the converged values of (\ref{consensus_yi}) and (\ref{consensus_zi}) be denoted as $r_i^\ast$ and $s_i^\ast$, respectively. Then, the optimal $\lambda^\ast$ in (\ref{lambda_optimal}) is computed as
\begin{align}
\lambda^\ast = -\frac{s_i^\ast}{r_i^\ast} \label{main2_lambda}
\end{align}
\end{corollary}
\begin{IEEEproof}
The proof is same as the proof of \textit{Theorem~\ref{theorem_energy_generation_main}}. 
\end{IEEEproof}
Now, after obtaining $\lambda^\ast$, we can accordingly compute $p_i^g$ as 
\begin{align}
p_i^g = -\frac{      \zeta_i + \lambda^\ast + {\lambda_i^v}^\ast}{2\xi_i}
\end{align}
Recursively, after obtaining $p_i^g$, we then replace $p_i^g$ in (\ref{lagrangian_lamv}) by new $p_i^g$ to repeat the overall process until $p_i^g$ converges to a value. The overall process for joint coordination is summarized in the \textit{Algorithm~\ref{jointcoordination}}.
\begin{algorithm}
\caption{Joint coordination $1$}\label{jointcoordination}
\begin{algorithmic}[1]
\BState \emph{Initial optimal generation}:
\State $\textit{select}$ $\epsilon_i$
\State $k \gets 0$
\State $\textit{obtain}~p_i^g(k)$ by solving (\ref{g_const1}) and (\ref{g_const2})
\Procedure{}{}
\State $i \gets \textit{patlen}$
\BState \emph{Coordination for energy flow}:
\State $\textit{given}~p_i^g(k)$, $\textit{obtain}~\lambda_i^v(k)$ by solving (\ref{lagrangian_pij}), (\ref{lagrangian_lamv}), (\ref{lagrangian_lame})
\BState \emph{Coordination for energy generation}:
\State $k \gets k+1$
\State $\textit{with~computed}~\lambda_i^v(k-1)$, $\textit{obtain}~p_i^g(k)$ by solving (\ref{g_const1}) and (\ref{additional})
\If {$\Vert p_i^g(k) - p_i^g(k-1) \Vert \leq \epsilon_i,~\forall i$,}
\State \textbf{goto} \emph{end}.
\State \textbf{close};
\EndIf
\State \textbf{otherwise}  $\textit{repeat}$~$\textbf{procedure}$      
\State
\EndProcedure
\State \textbf{end} \emph{stop}.
\end{algorithmic}
\end{algorithm}
However, although this approach seems efficient and can be implemented in a fully distributed way, the convergence of the recursive process is hard to prove. 

Alternatively, the energy generation and energy flow optimization algorithms might be solved simultaneously in a coupled manner. By inserting $\lambda$ of (\ref{lambda_optimal}) into (\ref{additional}), we are able to obtain $p_i^g$ that can be expressed as a function of $\lambda_i^v$. Then, we insert $p_i^g$ into (\ref{lagrangian_lamv}) to, eventually, have
\begin{align}
\sum_{j \in \mathbf{N}_i} \frac{1}{4 \alpha_{ij}} (\lambda_i^v -\lambda_j^v) = -\sum_{j \in \mathbf{N}_i} \left\{ \frac{\beta_{ij}}{2 \alpha_{ij}}      \right\} - \frac{\zeta_i}{2 \xi_i} \nonumber\\+ \frac{  \sum_{j=1}^n \left\{    \frac{\zeta_j}{2 \xi_j} +  \frac{\lambda_j^v}{2 \xi_j}   + p_j^d   \right\}         }{   \xi_i    \sum_{j=1}^n     \frac{1}{\xi_j} }  
- \frac{\lambda_i^v}{2 \xi_i}  - p_i^d \label{rhs_joint}
\end{align}
Using the substitution $\frac{\lambda_i^v}{2 \xi_i} = \frac{\lambda_i^v        \sum_{j=1}^n     \frac{1}{2 \xi_j}      }{2 \xi_i  \sum_{j=1}^n     \frac{1}{2 \xi_j}}$,
we can change the right-hand side of (\ref{rhs_joint}) as $
-\sum_{j \in \mathbf{N}_i} \left\{ \frac{\beta_{ij}}{2 \alpha_{ij}} \right\} - \frac{\zeta_i}{2 \xi_i}  -p_i^d 
 +\frac{  \sum_{j=1}^n \left\{    \frac{\zeta_j}{2 \xi_j}   + p_j^d   \right\}         }{   \xi_i    \sum_{j=1}^n \frac{1}{\xi_j} }  +\frac{  \sum_{j=1}^n \left\{  (\lambda_j^v - \lambda_i^v)/2\xi_j   \right\}         }{   \xi_i    \sum_{j=1}^n \frac{1}{\xi_j} }$.
Now, after defining adjacency matrix $A_{\mathcal{G}}^k$ for complete graph with weightings $\frac{1/(2 \xi_j)}{\xi_i    \sum_{j=1}^n \frac{1}{\xi_j}}$ for $(i,j)$ element,  we can make the Laplacian matrix $\mathcal{L}_{\mathcal{G}}^k$ for the complete graph. Then, (\ref{rhs_joint}) can be rewritten as
\begin{align}
(\mathcal{L}_{\mathcal{G}}+ \mathcal{L}_{\mathcal{G}}^k) \lambda^v = \mathbf{w}_J \label{joint_final_eq}
\end{align}
where $\mathbf{w}_J$ is the length $n$ vector with $i$-th element given as $-\sum_{j \in \mathbf{N}_i} \left\{ \frac{\beta_{ij}}{2 \alpha_{ij}} \right\} - \frac{\zeta_i}{2 \xi_i}  -p_i^d 
 +\frac{  \sum_{j=1}^n \left\{    \frac{\zeta_j}{2 \xi_j}   + p_j^d   \right\}         }{   \xi_i    \sum_{j=1}^n \frac{1}{\xi_j} }$. Thus, if $\mathcal{L}_{\mathcal{G}}+ \mathcal{L}_{\mathcal{G}}^k$ is a nonsingular matrix,  we can obtain $\lambda^v$ as $
 \lambda^v = (\mathcal{L}_{\mathcal{G}}+ \mathcal{L}_{\mathcal{G}}^k)^{-1} \mathbf{w}_J$.
However, note that the summation of two Laplacian matrices does not imply a positive definiteness; so there may not exist an inverse of  the matrix $\mathcal{L}_{\mathcal{G}}+ \mathcal{L}_{\mathcal{G}}^k$. Moreover, obtaining the inverse of $\mathcal{L}_{\mathcal{G}}+ \mathcal{L}_{\mathcal{G}}^k$ is a centralized approach and demands a huge computation. Thus, it is highly preferred to solve (\ref{joint_final_eq}) in a distributed way. To solve (\ref{joint_final_eq}) in a distributed way, we propagate it such as
\begin{align}
\dot{\lambda}^v = -(\mathcal{L}_{\mathcal{G}}+ \mathcal{L}_{\mathcal{G}}^k) \lambda^v + {\bf w}_J \label{laplacian_lambda_v_prop4}
\end{align}
Then, similarly to \textit{Theorem~\ref{theorem_flow}}, we could make the following theorem:
\begin{theorem} In (\ref{laplacian_lambda_v_prop4}), $\lambda^v$ converges to a constant vector, and the converged solution is equivalent to the solution of $(\mathcal{L}_{\mathcal{G}}+ \mathcal{L}_{\mathcal{G}}^k) \lambda^v - {\bf w}_J =0$. \label{theorem_joint_main}
\end{theorem}
\begin{IEEEproof}
Using $\mathcal{Q}$ and the vector $\overline{\lambda}^v$ defined in (\ref{Trans_Q}) and (\ref{new_lambda_trans_Q}), 
we transform (\ref{laplacian_lambda_v_prop4}) as
\begin{align}
\dot{\overline{\lambda}^v} = - \mathcal{Q}^T (\mathcal{L}_{\mathcal{G}}+ \mathcal{L}_{\mathcal{G}}^k) \mathcal{Q} \overline{\lambda}^v + \mathcal{Q}^T {\bf w}_J \label{laplacian_lambda_v_prop5}
\end{align}
It is clear that $ \sum_{j \in \mathbf{N}_1} \left\{ \frac{\beta_{1j}}{2 \alpha_{1j}} \right\}+\sum_{j \in \mathbf{N}_2} \left\{ \frac{\beta_{2j}}{2 \alpha_{2j}} \right\}+  \cdots         +\sum_{j \in \mathbf{N}_n} \left\{ \frac{\beta_{nj}}{2 \alpha_{nj}} \right\} =0$. Also, with some manipulations, we can get 
\begin{align}
\sum_{i=1}^n \left\{-\frac{\zeta_i}{2 \xi_i}  -p_i^d  +\frac{  \sum_{j=1}^n \left\{    \frac{\zeta_j}{2 \xi_j}   + p_j^d   \right\}         }{   \xi_i    \sum_{j=1}^n \frac{1}{\xi_j} } \right\} = 0
\end{align}
Thus,  we also see that $({\bf 1}_n)^T {\bf w}_J=0$, which means that we can change (\ref{laplacian_lambda_v_prop5}) as
\begin{equation}
\dot{\overline{\lambda}^v} = -\left[
\begin{array}{c|c}
  0 & 0 \\
  \hline
  0 & \mathcal{Q'}^T (\mathcal{L}_{\mathcal{G}}+ \mathcal{L}_{\mathcal{G}}^k) \mathcal{Q'}
\end{array}
\right]  \overline{\lambda}^v   +    \left[
\begin{array}{c}
  0  \\
  \hline
  {\bf w}_J' 
\end{array}
\right]       
\end{equation}
where $\mathcal{Q}' \in \Bbb{R}^{n \times (n-1)}$ is same as (\ref{q_prime}). Since the matrix $\mathcal{Q'}^T (\mathcal{L}_{\mathcal{G}}+ \mathcal{L}_{\mathcal{G}}^k) \mathcal{Q'}$ is positive definite, it is clear that ${\overline{\lambda}^v}$ is convergent. Consequently, by the \textit{Lemma~\ref{lemma_const}}, the proof is completed. 
\end{IEEEproof}
 
Note that the term $\mathcal{L}_{\mathcal{G}}^k \lambda^v$ in the right-hand side requires information of a complete graph; so it needs information from all nodes. To resolve this problem, let us suppose that $\dot{\lambda}^v$ in (\ref{laplacian_lambda_v_prop4}) is the derivative of ${\lambda}^v$ with respect to the time variable $t$. The sampling interval between two consecutive instants $t_1$ and $t_2$ is denoted as $\Delta t$. Let us choose another variable $\kappa$ and its sampling interval between two consecutive instants is denoted as $\Delta \kappa$. It is supposed that $\Delta \kappa \ll \Delta t$. 
Let us denote the third term in the right-hand side of (\ref{rhs_joint}) 
at time instant $t$ as
\begin{align}
\phi_i(t) \triangleq \frac{  \sum_{j=1}^n \left\{    \frac{\zeta_j}{2 \xi_j} +  \frac{\lambda_j^v}{2 \xi_j}   + p_j^d   \right\}         }{   \xi_i    \sum_{j=1}^n     \frac{1}{\xi_j} }
\end{align}
To find $\phi_i(t)$ in a distributed way at instant $t$, we use subsidiary variables $r_i^{\phi}(\kappa)$ and $s_i^{\phi}(\kappa)$ and assign their initial values, in terms of $\kappa$, as $r_i^{\phi}(0) = \frac{1}{\xi_i}        
$ and $s_i^{\phi}(0) =   \frac{\zeta_i}{2 \xi_i} +  \frac{\lambda_i^v}{2 \xi_i}   + p_i^d$ respectively. 
Let us update $r_i^\phi(\kappa)$ and $s_i^\phi(\kappa)$ by the following distributed consensus algorithms $
\frac{{\text d}{r}_i^\phi}{{\text d}\kappa} = -\sum_{j \in \mathbf{N}_i} (r_i^\phi(\kappa) - r_j^\phi(\kappa))$ and $\frac{{\text d}{s}_i^\phi}{{\text d}\kappa} = -\sum_{j \in \mathbf{N}_i} (s_i^\phi(\kappa) - s_j^\phi(\kappa))$.
Since we have selected $\Delta \kappa \ll \Delta t$, it can be supposed that the above consensus equations converge very quickly within $\Delta t$. Denoting the converged values at time instant $t$ as $r_i^{\phi\ast}(\kappa)$ and $s_i^{\phi\ast}(\kappa)$, we can obtain $\phi_i(t)$ as $
\phi_i^\ast(t) = \frac{s_i^{\phi\ast}(\kappa)}{\xi_i r_i^{\phi\ast}(\kappa)} \label{phi_i_t}$. Thus, $\lambda_i^v$ is propagated as 
\begin{align}
\dot{\lambda_i^v} &= - \sum_{j \in \mathbf{N}_i} \frac{1}{4 \alpha_{ij}} (\lambda_i^v -\lambda_j^v)  -\sum_{j \in \mathbf{N}_i} \left\{ \frac{\beta_{ij}}{2 \alpha_{ij}}      \right\} - \frac{\zeta_i}{2 \xi_i} \nonumber\\&~+  \phi_i^\ast(t)  
- \frac{\lambda_i^v}{2 \xi_i}  - p_i^d \label{lambda_joint_distributed_optimal}
\end{align}
which will be convergent according to the \textit{Theorem~\ref{theorem_joint_main}}.
The overall idea is summarized in the \textit{Algorithm~\ref{jointcoordination2}}.
\begin{algorithm}
\caption{Joint coordination $2$}\label{jointcoordination2}
\begin{algorithmic}[1]
\BState \emph{Initial optimal generation}:
\State $\textit{two variables}$ $t, \kappa$ ($\Delta \kappa \ll \Delta t$)
\State {at every} $t$, \textit{obtain} $\phi_i^\ast(t)$ in a distributed way 
\State $\textit{propagate}$  $\lambda_i^v$ according to  (\ref{lambda_joint_distributed_optimal})
\State with converged $\lambda_i^v$, \textit{compute} $\lambda$ in a distributed way 
using (\ref{lambda_optimal})
\State $\textit{obtain}~p_i^g$ by $p_i^g = -\frac{1}{2 \xi_i}[ \zeta_i + \lambda + \lambda_i^v]$
\State with generated $p_i^g$ and converged $\lambda_i^v$, $\textit{compute}$ $p_{ij}$ accordingly
\end{algorithmic}
\end{algorithm}
\begin{remark}
The joint coordination algorithm outlined in the \textit{Algorithm~\ref{jointcoordination2}} can be implemented in a fully distributed way; however, to realize this algorithm, we need two loops, i.e., inner loop operated by the variable $\kappa$ and outer loop operated by the variable $t$. To ensure convergence, we need to manage two loops such that the inner loop works much faster than the outer loop.   
\end{remark}

\section{Illustration By Simulation} \label{sec_simulation}
\begin{figure}
\includegraphics[width=0.45\textwidth]{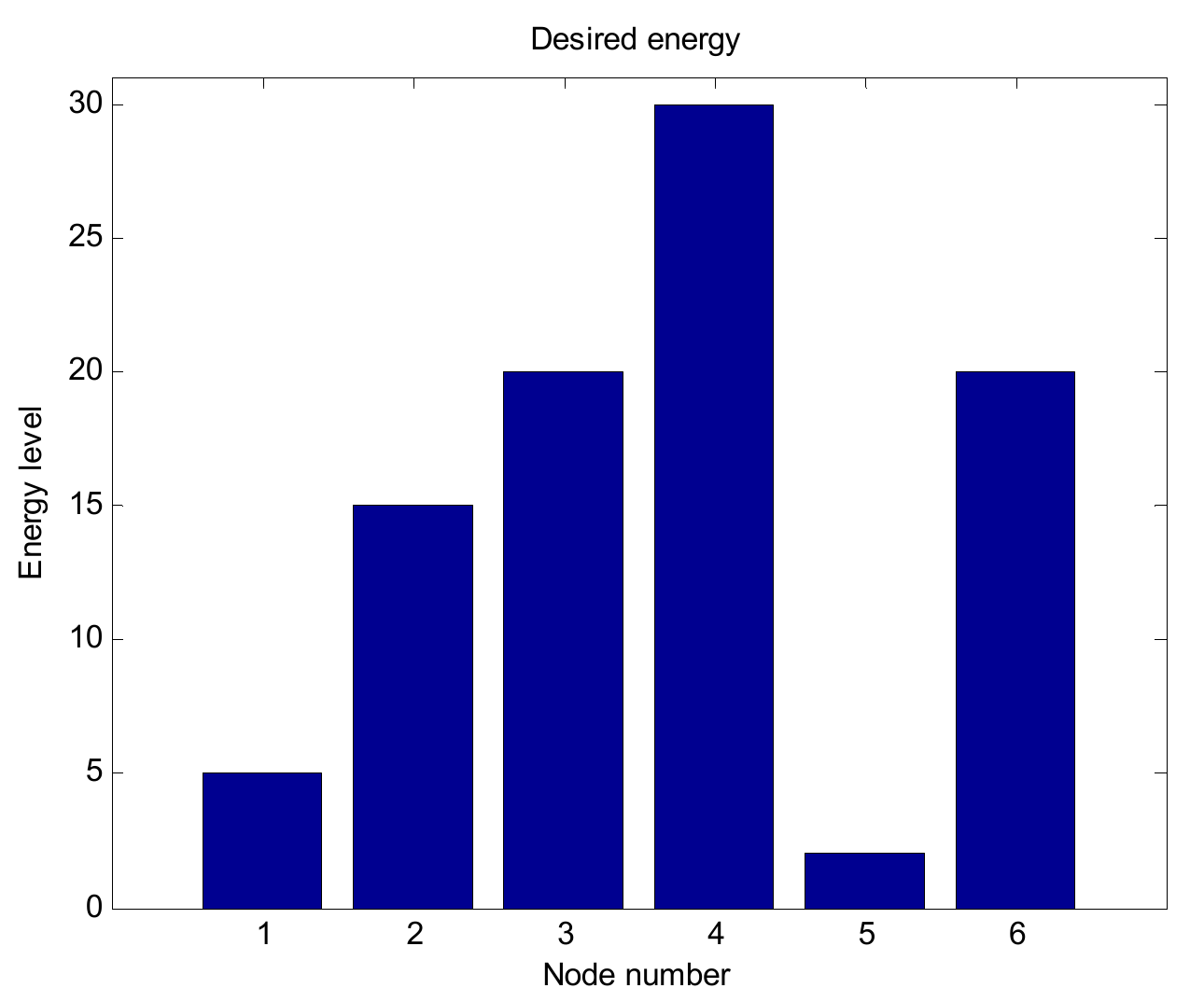}
\caption{Desired energy levels of each node.}
\label{fig_desired}
\end{figure}

For the verification of the algorithms, we conduct numerical simulations. We consider 6 nodes with the following edge set (the underlying graph is represented by Fig.~\ref{sgn}):
\begin{align}
\mathcal{E} \triangleq \{ (1,2), (2,3), (3,4), (1,5), (3,5), (4,5), (4,6) \}
\end{align}
The cost for each edge is assigned with the weights $\alpha_{1,2} =5, \alpha_{2,3} = 7, \alpha_{3,4}=3, \alpha_{1,5}=4, \alpha_{3,5}=6, \alpha_{4,5}=8, \alpha_{4,6}=7$ and $\beta_{1,2}=0.01, \beta_{2,3}=0.01, \beta_{3,4}=0.01, \beta_{1,5}=0.01, \beta_{3,5}=0.01, \beta_{4,5}=0.01, \beta_{4,6}=0.01$, and $\gamma_{i,j}=0,~\forall (i,j) \in \mathcal{E}$. For the selection of cost function for energy generation, we need to choose target parabola equation as $f(p_{i})= \xi_{i} p_{i}^2 + \eta_{i} p_{i} + \zeta_{i}$. Mainly, we may be interested in the location of the minimum point, and the slope of parabola equation. We can change $f(p_{i})$ as
\begin{align}
f(p_{i}) = \xi_{i} \left(p_{i} + \frac{\eta_{i}}{2\xi_{i}}\right)^2 -\frac{(\eta_{i})^2}{4 \xi_{i}} + \zeta_{i}
\end{align}
The location of the minimum point is $\left(-\frac{\eta_{i}}{2\xi_{i}}, -\frac{(\eta_{i})^2}{4 \xi_{i}} + \zeta_{i}   \right)$ and the slope is decided by $\xi_{i}$. For the simulation, we select $\xi_1=10, \xi_2=15, \xi_3=12, \xi_4=10, \xi_5=10, \xi_6=15$ and $\frac{\eta_{1}}{2\xi_{1}} =10, \frac{\eta_{2}}{2\xi_{2}} =15, \frac{\eta_{3}}{2\xi_{3}} =12, \frac{\eta_{4}}{2\xi_{4}} =10, \frac{\eta_{5}}{2\xi_{5}} =10, \frac{\eta_{6}}{2\xi_{6}} =15$, and we select $\zeta_i= \frac{(\eta_{i})^2}{4 \xi_{i}}$ to make the minimum be equal to zero. The initial energy of each node is assumed zero, and the diagram in Fig.~\ref{fig_desired} shows the desired energy level of each node. 

\begin{figure}
\includegraphics[width=0.45\textwidth]{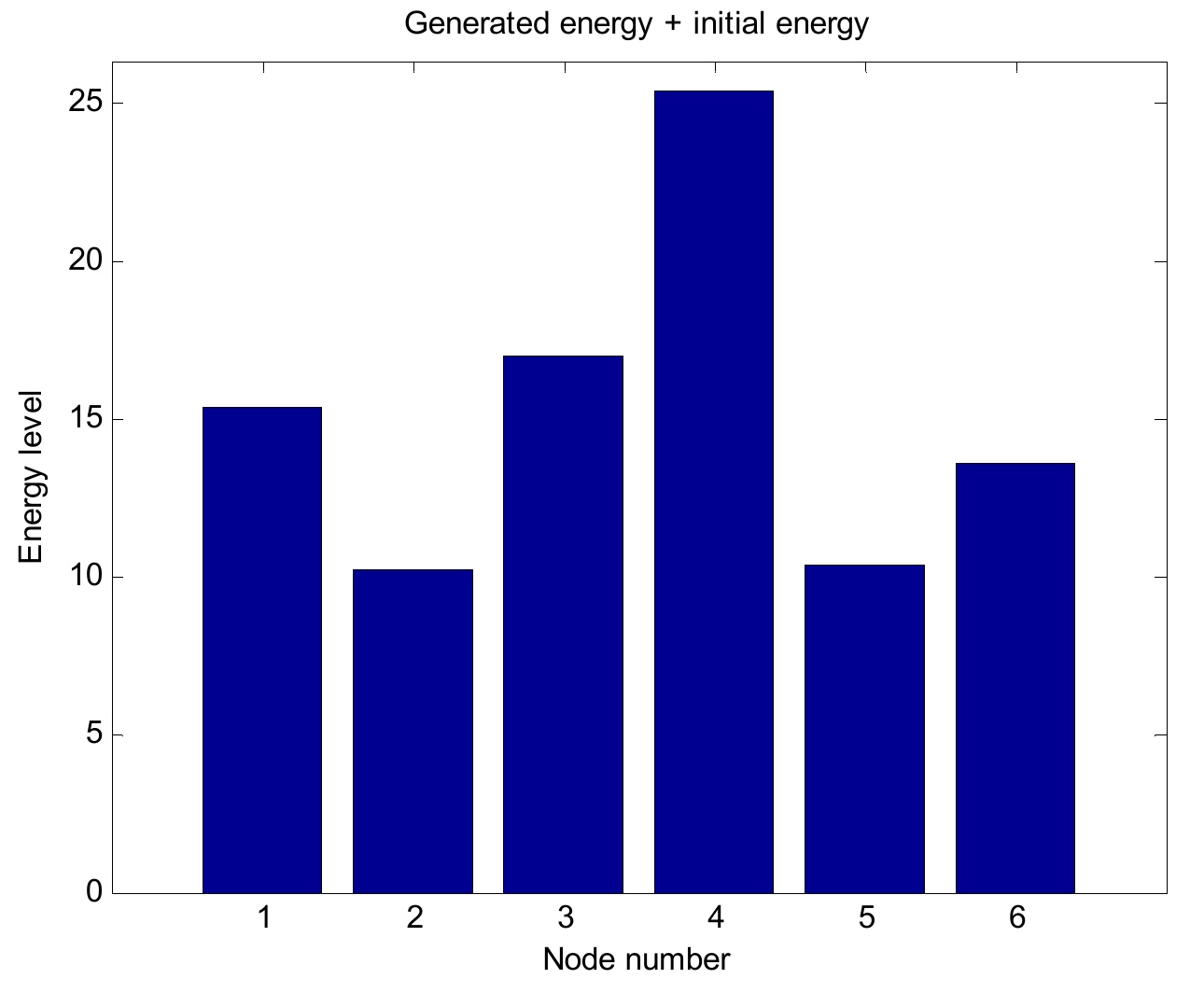}
\caption{Generated energy of each node via coordination of energy generation.}
\label{fig_generated1}
\end{figure}

\begin{figure}
\includegraphics[width=0.45\textwidth]{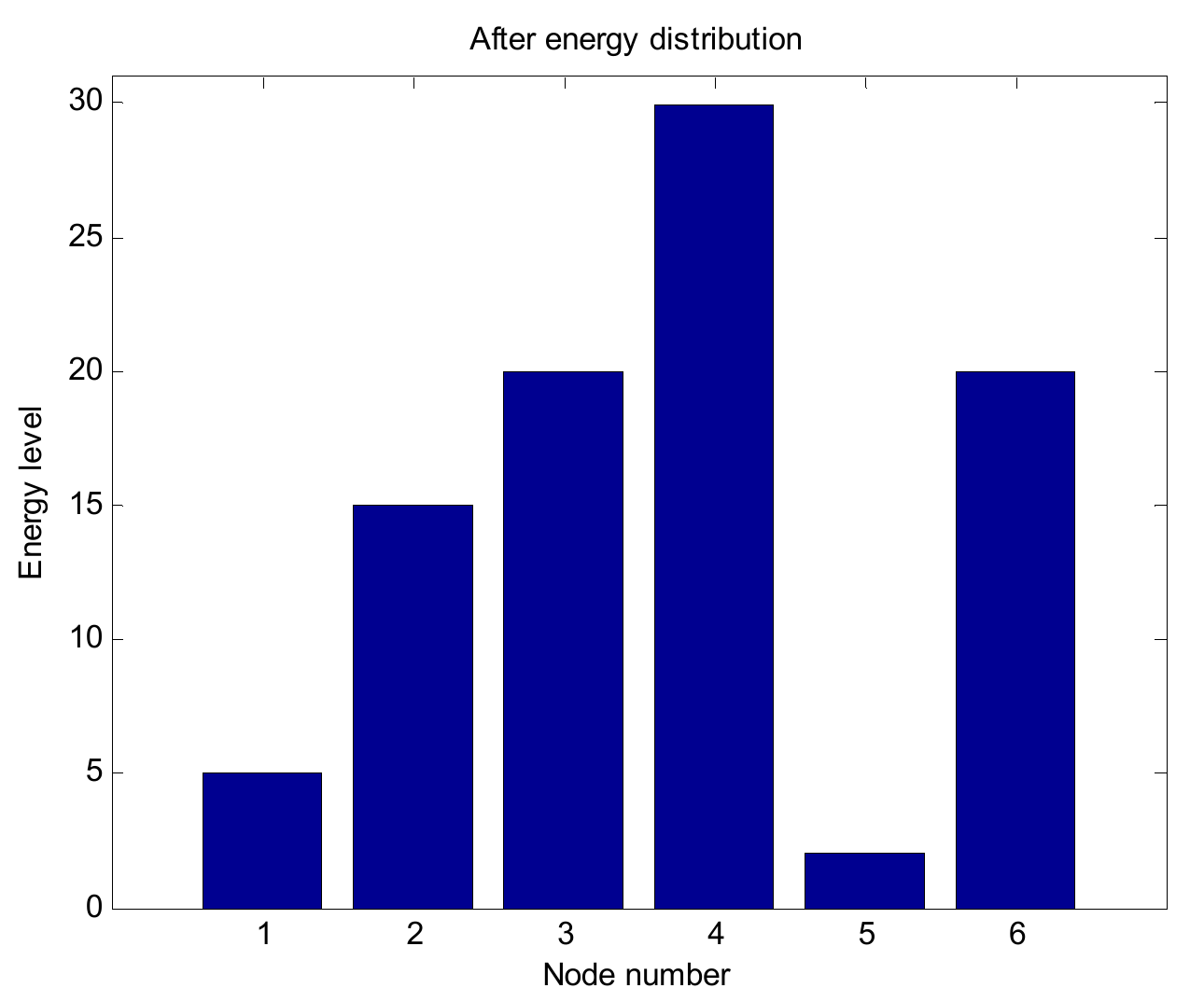}
\caption{Energy levels of each node after coordination.}
\label{fig_after}
\end{figure}

\begin{figure}
\includegraphics[width=0.45\textwidth]{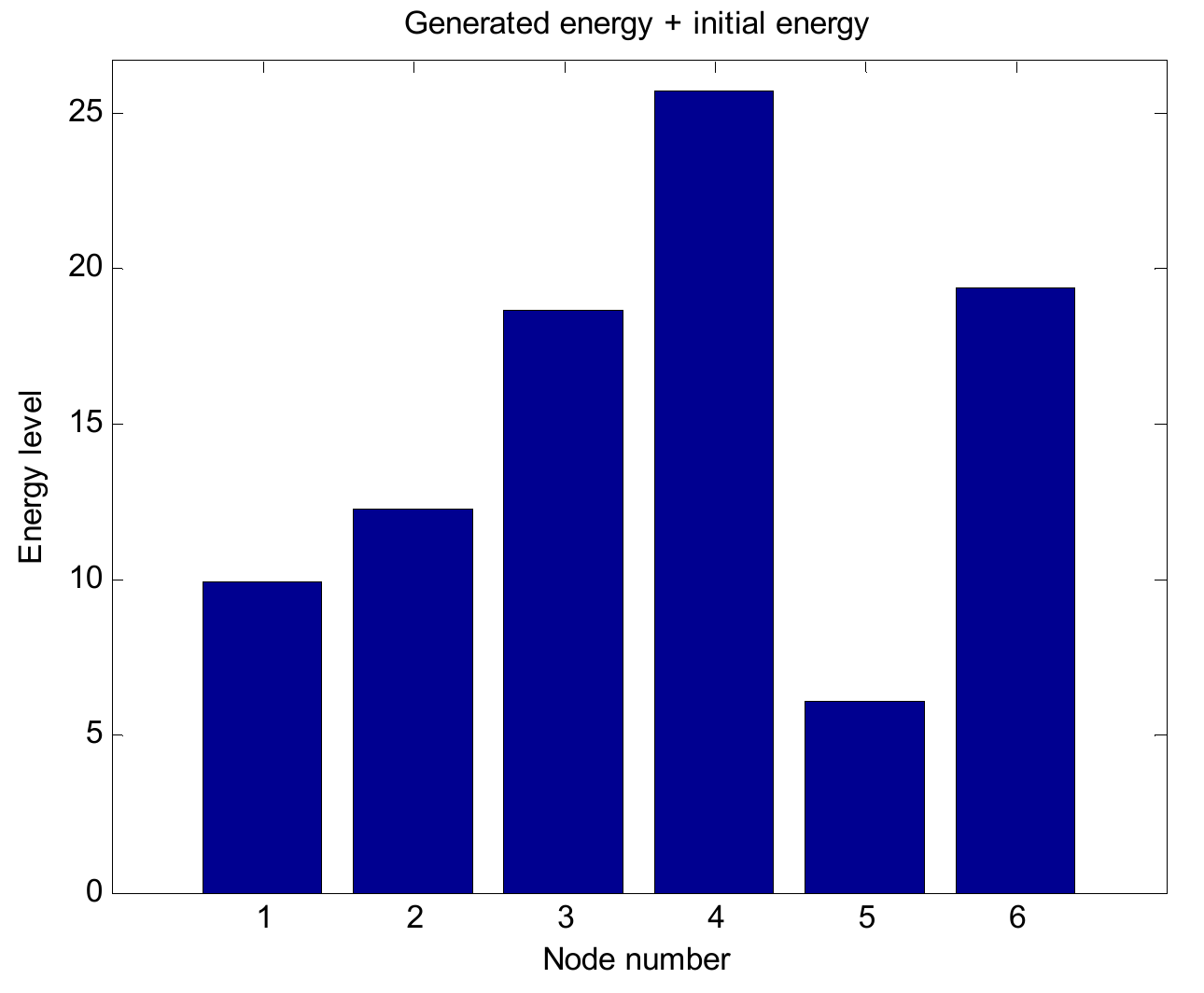}
\caption{Generated energy of each node via joint coordination.}
\label{fig_generated2}
\end{figure}

Using the coordination algorithm given in Subsection~\ref{subsec_gen}, we can generate the energy in each node in the direction of minimizing the cost function. Fig.~\ref{fig_generated1} shows the generated energy in each node after coordinated by the \textit{Algorithm~\ref{coordination1}} given in Subsection~\ref{subsec_gen}. After generating energy, we use the energy distribution algorithm (i.e., \textit{Algorithm~\ref{coordination2}} ) given in Subsection~\ref{subsec_flow}. Fig.~\ref{fig_after} shows the energy level after energy flow. As shown in this figure, the desired energy has been achieved at all nodes. It is noticeable that we had the exactly same results as Fig.~\ref{fig_generated1} and Fig.~\ref{fig_after} by using centralized optimization scheme (i.e., by directly solving the optimization problems using a centralized computation). Thus, we could see that the proposed algorithms solve the global optimization problems in a distributed way.  We also use the \textit{Algorithm~\ref{jointcoordination2}} for the joint coordination. Fig.~\ref{fig_generated2} shows the energy level generated by the joint coordination. Note that the generated energies by the joint coordination are different from the energy generated by the  coordination algorithm (i.e., \textit{Algorithm~\ref{coordination1}}) given in Subsection~\ref{subsec_gen}. 

The amount of energy generated by the \textit{Algorithm~\ref{coordination1}} is $p^g=(13.2903, 12.1935,$ $17.7419, 23.2903, 8.2903,   17.1935)$, and the amount of energy flows by the \textit{Algorithm~\ref{coordination2}} is $p_{12}=-4.9972$, $p_{15}=-3.2888$, $p_{23}=-2.1840$, $p_{34}=-4.4449$,  $p_{35}=4.5182$, $p_{45}=5.0548$, and $p_{46}=-2.7919$. The energy generated by the joint coordination is $p^g=(9.9318, 12.2600, 18.6629, 25.6667, 6.1108,$ $19.3678)$, and amount of flows is  $p_{12}=-3.4592$, $p_{15}=-1.4750$, $p_{23}=-0.7190$, $p_{34}=-2.1203$,  $p_{35}=2.7364$, $p_{45}=2.8468$, and $p_{46}=-0.6328$.

To check optimality, we compare the cost of the coordination for energy generation by  \textit{Algorithm~\ref{coordination1}} (i.e., $\sum_{i \in \mathcal{V}} (\xi_i (p_i^g)^2 + \zeta_i p_i^g + \eta_i)$, the cost of coordination for energy flow by  \textit{Algorithm~\ref{coordination2}} (i.e., $\sum_{(i,j) \in \mathcal{E}} (\alpha_{ij} p_{ij}^2 + \beta_{ij} p_{ij} + \gamma_{ij})$), and the cost of joint coordination by  
\textit{Algorithm~\ref{jointcoordination2}} (i.e., $\sum_{i \in \mathcal{V}} (\xi_i (p_i^g)^2 + \zeta_i p_i^g + \eta_i) + \sum_{(i,j) \in \mathcal{E}} (\alpha_{ij} p_{ij}^2 + \beta_{ij} p_{ij} + \gamma_{ij})$). By the coordination of energy generation, we have the cost of $\sum_{i \in \mathcal{V}} (\xi_i (p_i^g)^2 + \zeta_i p_i^g + \eta_i)=559.3548$ for generation, and by the coordination of energy flow, the cost for flow is $\sum_{(i,j) \in \mathcal{E}} (\alpha_{ij} p_{ij}^2 + \beta_{ij} p_{ij} + \gamma_{ij})=1284.3$; so, the summation of the overall costs is $1843.7$ in decoupled coordination approach (i.e., by using the \textit{Algorithm~\ref{coordination1}} and \textit{Algorithm~\ref{coordination2}} in serial). By the joint coordination, the generation cost is $857.2791$ and the cost for flow is $396.3495$; so the total cost is $1253.6$ in jointly coupled coordination. From these results, we can see that even though the generation cost from the joint coordination is higher than the cost by the coordination of energy generation, the cost for flow is significantly reduced by the joint coordination. Thus, in this example, the overall cost has been reduced as much as 32$\%$ by the joint coordination. So, we can see that the joint coordination algorithm has reduced the cost arising from the flow by way of paying more cost in the energy generation process.

\section{Conclusion}            \label{sec_conclusion}
In this paper, we have proposed three coordination laws for energy generation and distribution in smart grid energy networks. Energy is coordinated to generate in each node and flow through physical layer with information obtained through communication layer. 
The coordination is conducted in each node only using local relative information. 
As the main claim of this paper, even though the coordination is conducted using local information, the optimal generation and distribution are achieved in a global sense. As illustrated through simulations, when the coordinations for energy generation and energy distribution are conducted in a decoupled way without taking account of the interactive characteristics between them, the coordination laws attempt to minimize the cost only in terms of generation or in terms of distribution. However, when the coordination is conducted in a coupled way considering the generation and distribution simultaneously, the overall cost is minimized globally. However, the two approaches have advantages and disadvantages respectively. When the coordinations are conducted in a decoupled way, each node needs less communication with its neighboring nodes; but when the coordination is conducted in a coupled way, each node needs more information exchanges with its neighboring nodes. Thus, depending upon application situations and/or main concerns,
we may have to select and use appropriate coordination law; for example, when energies have been generated already from distributed energy resources, we may want to use only the energy distribution law.    

\section*{Acknowledgement}
This research was supported by Korea Electric Power Corporation through Korea Electrical Engineering \& Science Research Institute (Grant number: R15XA03-47).

\bibliographystyle{IEEETran}
\bibliography{energy_distribution_bib} 
\end{document}